\documentclass[a4paper,12pt]{amsart}
\usepackage{tabularx}
\usepackage[dvips]{graphicx}
\usepackage{amsmath}
\usepackage{amssymb}
\usepackage{amsbsy}

\usepackage{amsgen}
\usepackage{amsopn}
\usepackage{amscd}
\usepackage{amstext}
\usepackage{amsthm}

\theoremstyle{definition}
\newtheorem{definition}{Definition}[section]
\theoremstyle{plain}
\newtheorem{theorem}{Theorem}[section]
\newtheorem{lemma}{Lemma}[section]

\newtheorem{proposition}{Proposition}[section]

\theoremstyle{remark}
\newtheorem{remark}{Remark}[section]

\theoremstyle{definition}

\newtheorem{acknow}{Acknowledgment}

\title{Homotopy classification of nanophrases \\
with less than or equal to four letters}

\author{FUKUNAGA Tomonori}
\begin{document}

\maketitle

\begin{abstract}
In this paper we give the stable classification of 
ordered, pointed, oriented multi-component curves on surfaces with
 minimal crossing number less than or equal to 2 such that any 
equivalent curve has no simply
 closed curves in its components. To do this, we use the
theory of words and phrases which was introduced by V. Turaev. Indeed we 
give the homotopy classification of nanophrases with less than or equal to
 4 letters. It is an extension of the classification of 
nanophrases of length 2 with less than or equal to 4 letters which 
was given by the author in a previous paper.     
\end{abstract}

{\bf keywords:} Nanophrases, Homotopy, Multi-component curves, Stable
equivalent \\ \par

Mathematics Subject Classification 2000: Primary 57M99; Secondary 68R15

\markboth{FUKUNAGA Tomonori}{Homotopy Classification of Nanophrases}

\section{Introduction.} 
The study of curves via words was introduced by C. F. Gauss \cite{ga}. 
Gauss encoded closed planar curves by words of 
certain type which are now called Gauss words.
We can apply this method to encode multi-component curves on surfaces.
For instance, in \cite{tu2} and \cite{tu3} V. Turaev studied 
stable equivalence classes of curves on surfaces by using 
generalized Gauss words (called nanowords). \par
More precisely a nanoword over an alphabet $\alpha$ 
endowed with an involution 
$\tau:\alpha \longrightarrow \alpha$ is a word in an alphabet 
$\mathcal{A}$ endowed with a projection $\mathcal{A} \ni A \mapsto
|A| \in  \alpha$ such that every letter appears twice or not at all.
In the case where the alphabet $\alpha$ consists of two elements 
permuted by $\tau$, the notion of a nanoword over $\alpha$ is equivalent
to the notion of an open virtual string introduced in \cite{tu5}.  
\par
Turaev introduced the homotopy equivalence
on the set of nanowords over
$\alpha$. The homotopy equivalence relation 
is generated by three types of moves on nanowords.
The first move consists of deleting two consecutive entries of the same
letter. The second move has the form $xAByBAz \mapsto xyz$ where $x,y,z$ 
are words and $A,B$ are letters such that $|A|=\tau(|B|)$.
The third move has the form $xAByACzBCt \mapsto xBAyCAzCBt$ where
$x,y,z,t$ are words and $A,B,C$ are letters such that $|A|=|B|=|C|$.
These moves are suggested by the three local deformations of 
curves on surfaces (See Fig. \ref{fig1} and \cite{tu2} for more
details).
In \cite{tu2} Turaev showed that 
a stable equivalence class of an oriented 
pointed curve on a surface is identified with 
a homotopy class of nanoword in a 2-letter alphabet. Moreover 
Turaev extended this result to multi-component curves. In fact a stable
equivalence class of an oriented, ordered, pointed 
multi-component curve on a surface is identified  with a
homotopy class of a nanophrase
in a 2-letter alphabet. Roughly speaking, a nanophrase is a sequence of 
words where concatenation of those words is a nanoword
(See also sub-section \ref{nphhom} and section \ref{nvm} for more
details). Thus, using Turaev's theory of words and phrases, we can 
treat curves on surfaces algebraically.\par
Homotopy classification of nanowords was given by Turaev in \cite{tu1}.
Turaev gave the classification of nanowords less than or equal to 
$6$ letters. Moreover, the author introduced new invariants of 
nanophrases and gave the  
homotopy classification of nanophrases of length $2$ with less than or 
equal to $4$ letters in \cite{fu}, using Turaev's classification of 
nanowords.\par 
The purpose of this paper is to  give the classification 
theorem of nanophrases over arbitrary alphabet with less than or equal
$4$ letters without the condition on length.
As a corollary of this theorem, we classify the multi-component curves
with minimum crossing number less than or equal to $2$ which has no
``untide'' components up to stable equivalence (Theorem \ref{mthm}).\par
The constitution of this paper is as follows. In sections 2-4 we
review the theory of multi-component curves and  
the homotopy theory of words and phrases. In section 5 we introduce 
known results on the classification of nanowords and nanophrases up to
homotopy and we generalize these results to phrases 
of an arbitrary length.
Finally in section 6 we give the proof of the main theorem 
in this paper.  
  
\section{Stable Equivalence of Multi-component Curves.}
\subsection{Multi-component curves.}
In this paper a \emph{curve} means the image of a 
generic immersion of an oriented
circle into an oriented surface. The word ``generic'' means that the
curve has only a finite set of self-intersections which are all double
 and transversal. A \emph{$k$-component curve} is defined in the same
 way as a curve with the difference that they may be formed by $k$
 curves rather than only one curve. These curves are \emph{components}
 of the $k$-component curve. A $k$-component curves are \emph{pointed}
if each component is endowed with a base point (the origin) distinct
from the crossing points of the $k$-component curve. A $k$-component
curve is \emph{ordered} if its components are numerated.Two ordered,
pointed curves are \emph{stably homeomorphic} if there is an orientation
 preserving homeomorphism of their regular neighborhoods in the ambient
surfaces mapping the first multi-component curve onto the second one and
preserving the order, the origins, and the orientations of the
components.\par
\begin{figure}
\centerline{\includegraphics[width=9cm]{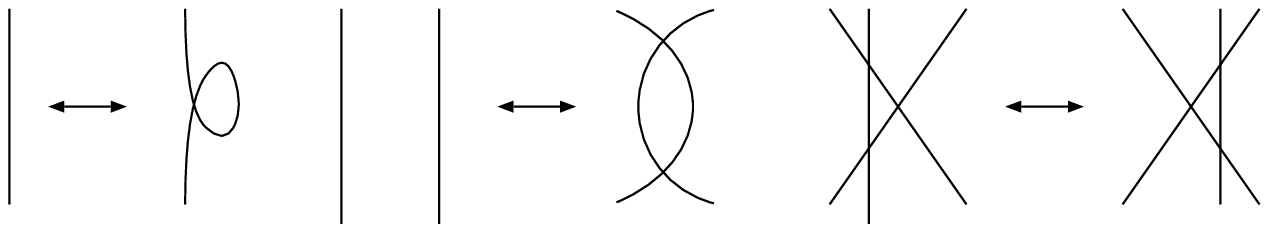}}
\caption{Three local deformations of curves.}\label{fig1}
\end{figure}

Now we define stable equivalence of ordered, pointed multi-component
curves \cite{ka2}: Two ordered, pointed multi-component 
curves are \emph{stably equivalent} if they can be related by a finite
sequence of the following transformations: (i) a move
replacing a ordered, pointed multi-component curve with a 
stably homeomorphic 
one; (ii) a deformation of a pointed curve in its ambient surface away
from the origin (such a deformation may push a branch of the
multi-component curves across another branch or a double point but not
across the origin of the curves) as in Fig. \ref{fig1}.

\par 
We denote the set of stable equivalence classes of ordered, pointed 
$k$-component curves by $\mathcal{C}_k$. 

\begin{remark}
The theory of stable equivalence class of multi-component curves on
 surfaces is closely related to the theory of virtual strings. See
\cite{kad} and \cite{tu5} for more details.
\end{remark}

We will show a following theorem by using Turaev's theory of words.

An ordered, pointed multi-component surface-curve is called
\emph{irreducible} if it is not stably equivalent to a surface-curve 
with a simply closed component. 

\begin{theorem}\label{mthm}
Any irreducible ordered, pointed multi-component surface-curve with 
minimal crossing number less than or equal to 2
is stably equivalent to one of the ordered, 
pointed multi-component curves arise from the following list (see
also Remark \ref{rem1}).
There are exactly 52 stable equivalence classes of irreducible 
ordered, pointed, multi-component surface-curves.
\begin{figure}
\centerline{\includegraphics[width=8cm]{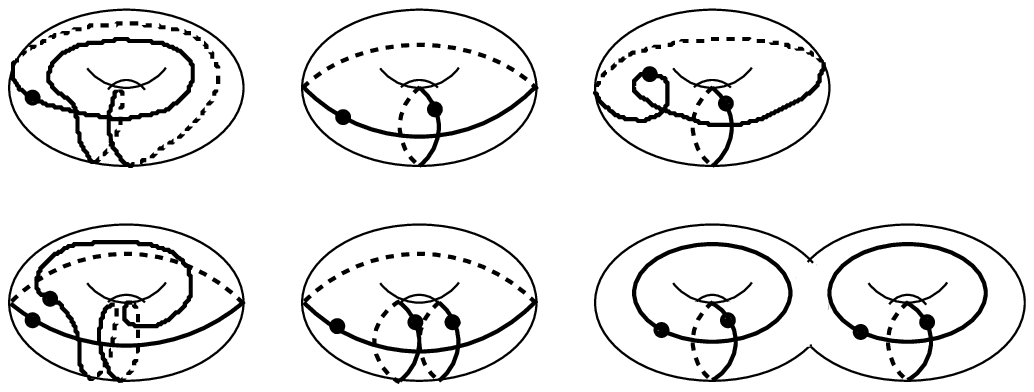}}
\caption{The list of curves.}\label{listfig}
\end{figure} 
\end{theorem}

\begin{remark}\label{rem1}
We want to list up the stable equivalence classes of irreducible
ordered, pointed multi-component surface-curves with 
minimal crossing number less than or equal to $2$. 
However there are too many curves to list up.
So in Fig. \ref{listfig} we make just the 
list of multi-component curves without order and orientation
 of the components. 
If we choose order and orientation, then 
we obtain a ordered, pointed multi-component curve.
Two different pictures from Fig. \ref{listfig} never produce
equivalent ordered, pointed multi-component surface-curves.
On the other hand it is possible that two different additional
structures (orientation and the order) on the same picture yield
equivalent ordered, pointed multi-component surface-curves.
More precisely, 2 (respectively 2, 8, 4, 24, 12) different ordered, pointed
multi-component surface-curves arise from the 
upper left (respectively upper middle, upper right, lower left, lower
middle, lower right) picture. 
By the Theorem \ref{ctk4}, ordered, pointed multi-component 
surface-curves arise from pictures in Fig. \ref{listfig} are stably equivalent 
if and only if nanophrases associated these curves are homotopic, and
we can obtain all of the stable equivalent classes of irreducible 
ordered, pointed
multi-component surface-curves with minimal crossing number 
less than or equal to $2$ by specifying order and orientation for
multi-component curves in Fig. \ref{listfig}.
\end{remark}

To prove the Theorem \ref{mthm}, we use Turaev's theory of words and phrases
which was introduced by V. Turaev in \cite{tu1} and \cite{tu2}.

\section{Turaev's Theory of Words and Phrases.}
In this section we review the theory of topology of words and phrases.
\subsection{Nanowords and their homotopy.}
An \emph{alphabet} is a set and \emph{letters} are its elements. 
A \emph{word of length $n \ge 1$ on an alphabet} $\mathcal{A}$ is a
mapping $w:\hat{n} \rightarrow \mathcal{A}$ where
$\hat{n}=\{1,2,\cdots,n\}$.
We denote a word of length $n$ by the sequence of letters $w(1)w(2)\cdots w(n)$.
A word $w:\hat{n} \rightarrow \mathcal{A}$ is a \emph{Gauss word} if
each element of $\mathcal{A}$ is the image of precisely two elements of
$\hat{n}$.
\par
For a set $\alpha$, an $\alpha$-alphabet is a set $\mathcal{A}$ 
endowed with a mapping $\mathcal{A} \rightarrow \alpha$ called
\emph{projection}.
The image of $A \in \mathcal{A}$ under this mapping is denoted $|A|$.
An \emph{\'etale word} over $\alpha$ is a pair 
(an $\alpha$-alphabet $\mathcal{A}$, a word on $\mathcal{A}$). 
A \emph{nanoword over} $\alpha$ is a pair
 (an $\alpha$-alphabet $\mathcal{A}$, a Gauss word on $\mathcal{A}$). 
We call an empty \'etale word in an empty $\alpha$-alphabet 
the \emph{empty nanoword}.
It is written $\emptyset$ and has length 0.
\par A \emph{morphism} of $\alpha$-alphabets
$\mathcal{A}_1$, $\mathcal{A}_2$ is a set-theoric mapping 
$f:\mathcal{A}_1 \rightarrow \mathcal{A}_2$ such that $|A|=|f(A)|$ 
for all $A \in \mathcal{A}_1$. If $f$ is bijective, then this morphism 
is an \emph{isomorphism}. Two \'etale words $(\mathcal{A}_1,w_1)$ and
$(\mathcal{A}_2,w_2)$ over $\alpha$ are \emph{isomorphic} if there is 
an isomorphism $f:\mathcal{A}_1 \rightarrow \mathcal{A}_2$ such that 
$w_2 = f \circ w_1$.
\par To define homotopy of nanowords we fix a finite set $\alpha$ 
with an involution $\tau : \alpha \rightarrow \alpha$ 
and a subset $S \subset \alpha \times \alpha \times \alpha$. 
We call the pair $(\alpha , S) $ \textit{homotopy data}.

\begin{definition}
 Let $(\alpha , S)$ be homotopy data. We define \textit{homotopy
 moves} (1) - (3) as follows:
 \par
 (1) $(\mathcal{A} , xAAy) \longrightarrow 
(\mathcal{A} \setminus \{ A \} , xy)$ 
 \par \enskip \enskip \enskip 
for all $A \in \mathcal{A}$ and $x,y$ are words in 
$\mathcal{A} \setminus \{ A \}$ such that $xy$ is a Gauss word.
 \par
(2) $(\mathcal{A} , xAByBAz) \longrightarrow (\mathcal{A} \setminus \{ A , B \} , xyz)$
   \par \enskip \enskip \enskip 
if $A , B \in \mathcal{A}$ satisfy $|B| = \tau (|A|)$. $x,y,z$ are words in 
$\mathcal{A} \setminus \{A,B\}$ such that \par \enskip \enskip \enskip 
 $xyz$ is a Gauss word.
 \par
(3) $(\mathcal{A} , xAByACzBCt) \longrightarrow (\mathcal{A} , xBAyCAzCBt)$
 \par \enskip \enskip \enskip if $A,B,C \in \mathcal{A}$ 
satisfy $(|A|,|B|,|C|) \in S $. $x,y,z,t$ are words in
 $\mathcal{A}$ such that \par \enskip \enskip \enskip $xyzt$ is a Gauss word.
\end{definition}

\begin{definition}
 Let $(\alpha , S)$ be homotopy data.
 Then nanowords $(\mathcal{A}_1 ,w_1)$ and $(\mathcal{A}_2 , w_2)$
 over $\alpha$ are \textit{$S$-homotopic}
 (denoted $(\mathcal{A}_1 , w_1) \simeq_{S} (\mathcal{A}_2 , w_2)$)
 if $(\mathcal{A}_2 , w_2)$ can be obtained from $(\mathcal{A}_1, w_1)$ by a
 finite sequence of isomorphism, $S$-homotopy moves 
 (1) - (3) and the inverse of moves (1) - (3).
\end{definition}

The set of $S$-homotopy classes of nanowords over $\alpha$ is denoted
$\mathcal{N}(\alpha,S)$.\par
To define $S$-homotopy of \'etale words we define
\emph{desingularization} of \'etale words $(\mathcal{A},w)$ over
$\alpha$
as follows: Set
$\mathcal{A}^d:=\{A_{i,j}:=(A,i,j) |A \in \mathcal{A}, 1 \le i < j \le
 m_w(A) \}$ with projection $|A_{i,j}|:=|A| \in \alpha$ for all
 $A_{i,j}$ (where $m_w(A) := Card(w^{-1}(A))$ ). 
The word $w^d$ is obtained from $w$ by first deleting all $A \in
 \mathcal{A}$ with $m_w(A)=1$. Then for each $A \in \mathcal{A}$ with
 $m_w(A)\ge 2$ and each $i = 1,2, \ldots m_w(A)$, we replace the $i$-th
 entry of $A$ in $w$ by $$A_{1,i}A_{2,i} \ldots
 A_{i-1,i}A_{i,i+1}A_{i,i+2} \ldots A_{i,m_w(A)}.$$
The resulting $(\mathcal{A}^d,w^d)$ is a nanoword of length 
$\sum_{A \in \mathcal{A}}
m_w(A)(m_w(A)-1)$ and called a \emph{desingularization of $(\mathcal{A},w)$}.
Then we define $S$-homotopy of \'etale words as follows:

\begin{definition}
 Let $w_1$ and $w_2$ be \'etale words over $\alpha$. Then $w_1$ and $w_2$
 are $S$-homotopic if $w_1^d$ and $w_2^d$ are $S$-homotopic.
\end{definition}

\subsection{Nanophrases and their homotopy.}\label{nphhom}  
In \cite{tu2}, Turaev proceeded similar arguments for phrases (sequence
of words).

\begin{definition}
A \emph{nanophrase} $(\mathcal{A},(w_1|w_2|\cdots|w_k))$ of 
length $k \ge 0$ over a set $\alpha$ is a 
pair consisting of an $\alpha$-alphabet $\mathcal{A}$ and a sequence of
 $k$ words $w_1,\cdots,w_k$ on $\mathcal{A}$ such that $w_1w_2\cdots
 w_k$ is a Gauss word on $\mathcal{A}$. We denote it simply by
 $(w_1|w_2|\cdots|w_k)$.
\end{definition}

\par By definition, there is a unique \emph{empty nanophrase} of length
$0$ (the corresponding $\alpha$-alphabet $\mathcal{A}$ is an empty set).

\begin{remark}
We can consider a nanoword $w$  
to be a nanophrase $(w)$ of length $1$.
\end{remark}

\par
A mapping $f:\mathcal{A}_1 \rightarrow \mathcal{A}_2$ is 
\emph{isomorphism} of two nanophrases if $f$ is an isomorphism of
$\alpha$-alphabets transforming the first nanophrase into the second one.
\par
Given homotopy data $(\alpha,S)$, we define homotopy moves on 
nanophrases as in section 3.1 with the only difference that the
2-letter subwords $AA$, $AB$, $BA$, $AC$ and $BC$ modified by these
 moves may occur in different words of phrase. 
Isomorphism and homotopy moves generate 
an equivalence relation $\simeq_S$ of $S$-homotopy on the classes of 
nanophrases over $\alpha$. We denote the set of $S$-homotopy classes of
nanophrases of length $k$ by $\mathcal{P}_k(\alpha,S)$.\par

\section{Nanophrases versus Multi-component Curves}\label{nvm}
In \cite{tu2}, Turaev showed that the special case of the study of homotopy
theory of nanophrases is equivalent to the study of $\mathcal{C}_k$.
More precisely, Turaev showed following theorem.

\begin{theorem}\emph{(Turaev \cite{tu2}).}
Let $\alpha_0$ is the set $\{a,b\}$ with involution $\tau:\alpha_0
 \rightarrow \alpha_0$ permuting $a$ and $b$, and $S_0$ is the diagonal of 
$\alpha_0 \times \alpha_0 \times \alpha_0$. Then
there is a canonical bijection $\mathcal{C}_k$ to 
$\mathcal{P}_k(\alpha_0,S_0)$. 
\end{theorem}

The method of making nanophrase $P(C)$ 
from ordered, pointed $k$-component curve $C$ is as
follows.
Let  us label the double points of $C$ by distinct letters
$A_1,\cdots,A_n$. Starting at the origin of first component of $C$ and
following along $C$ in the positive direction, we write down the labels 
of double points which we passes until the return to the origin.
Then we obtain a word $w_1$.
Similarly we obtain words $w_2,\cdots,w_k$ 
on the alphabet $\mathcal{A}=\{A_1,\cdots,A_n\}$ from 
second component, $\cdots$, $k$-th component.
Let $t_i^1$ (respectively, $t_i^2$) be the tangent vector 
to $C$ at the double point
labeled $A_i$ appearing at the first (respectively, second) passage
through this point. Set $|A_i|=a$, if the pair $(t_i^1,t_i^2)$ is
positively oriented, and $|A_i|=b$ otherwise. Then we obtain a required
 nanophrase $P(C):=(\mathcal{A},(w_1|\cdots|w_k))$.
 
By the above theorem if we classify the homotopy classes of
nanophrases, then we obtain the classification of ordered, pointed
multi-component curves under the stable equivalence as a corollary.

\begin{remark}
In \cite{SW}, D. S. Silver and S. G. Williams studied open virtual 
multi-strings. The theory of open virtual multi-strings is equivalent to 
the theory of pointed multi-component surface-curves. Silver and Williams
constructed invariants of open virtual multi-strings.   
\end{remark}

\section{Classification of Nanophrases.}
In this section, we give the homotopy classification of nanophrases with
less than or equal to $4$ letters under the assumption that a homotopy
data $S$ is the diagonal. In the remaining part of the paper we always
assume that homotopy data is the diagonal. Note that this assumption is not 
obstruct the our purpose.\par

\subsection{The case of  nanophrases of length $1$.}
The case of nanophrases of length $1$ (in other words the case
 of nanowords), Turaev gave the following classification theorem.

\begin{theorem}\label{tu}\emph{(Turaev \cite{tu1}).}
Let $w$ be a nanoword of length $4$ over $\alpha$.
 Then $w$ is either homotopic to the empty nanoword or isomorphic to the 
nanoword $w_{a,b}:=(\mathcal{A}=\{A,B\},ABAB)$ where 
 $|A|=a,|B|=b \in \alpha$ with $a \neq \tau(b)$.
 Moreover for $a \neq \tau(b)$, the nanoword $w_{a,b}$ is
 non-contractible and two nanowords $w_{a,b}$ and
 $w_{a^\prime,b^\prime}$ are homotopic if and only if 
$a=a^\prime$ and $b=b^\prime$.  
\end{theorem}    
 
\begin{remark}
In the paper \cite{tu1}, Turaev gave the classification of nanowords of
 length $6$. But in this paper we do not use this result. 
Classification problem of nanowords of length more than or equal
to $8$ is still open (See \cite{tu3}) .   
\end{remark} 
    
\subsection{The case of nanophrases of length $2$.} 
First we prepare following notations: $P_a:=(A|A)$, 
$P_{a,b}^{4,0}:=(ABAB|\emptyset),$
$P_{a,b}^{3,1}:=(ABA|B),$
$P_{a,b}^{2,2I}:=(AB|AB),$  
$P_{a,b}^{2,2II}:=(AB|BA),$
$P_{a,b}^{1,3}:=(A|BAB)$ and
$P_{a,b}^{0,4}:=(\emptyset|ABAB)$ 
with $|A|=a, \ |B|=b \in \alpha$.
If $a=\tau(b)$, then 
$P_{a,b}^{4,0}$, $P_{a,b}^{2,2I}$, $P_{a,b}^{2,2II}$ and
$P_{a,b}^{0,4}$ are homotopic to $(\emptyset|\emptyset)$. So in this paper, 
if we write $P_{a,b}^{4,0}$, 
$P_{a,b}^{2,2I}$, $P_{a,b}^{2,2II}$, $P_{a,b}^{0,4}$ then we always 
assume that $a \neq \tau (b)$. \par
In \cite{fu}, the author gave the classification of nanophrases of 
length $2$ with less than or equal to $4$ letters.

\begin{theorem}
Let $P$ be a nanophrase of length 2 with 2 letters. Then 
$P$ is not homotopic to $(\emptyset|\emptyset)$ if and only if 
$P$ is isomorphic to $P_a$.
Moreover $P_a$ and $P_{a^\prime}$ are homotopic if and only if $a=a^\prime$.
\end{theorem}

\begin{theorem}\label{ct24}
Let $P$ be a nanophrase of length 2 with 4 letters, then $P$ is
homotopic to $(\emptyset|\emptyset)$ or homotopic to
nanophrases of length 2 with 2 letters 
or isomorphic to one of the following nanophrases: $P_{a,b}^{4,0}$, $P_{a,b}^{3,1}$, 
$P_{a,b}^{2,2I}$, $P_{a,b}^{2,2II}$, $P_{a,b}^{1,3}$, $P_{a,b}^{0,4}$. 
For $(i,j) \in \{(4,0),(3,1),(2,2I),(2,2II),(1,3), (0,4)\}$ and any 
$a,b \in \alpha$, the nanophrase $P_{a,b}^{i,j}$ is neither homotopic to
$(\emptyset|\emptyset)$ nor homotopic to 
nanophrases of length 2 with 2 letters. The nanophrases $P_{a,b}^{i,j}$
 and $P_{a^\prime,b^\prime}^{i,j}$ are homotopic if and only if 
$a=a^\prime$ and $b=b^\prime$. For $(i,j) \neq (i^\prime,j^\prime),$ the 
nanophrases $P_{a,b}^{i,j}$ and
 $P_{a^\prime,b^\prime}^{i^\prime,j^\prime}$ are not homotopic for any 
$a,b,a^\prime,b^\prime \in \alpha$.
\end{theorem}

In this paper, we give the classification of nanophrases of length more 
than or equal to $3$ with $4$ letters.

\subsection{Homotopy invariants of nanophrases.}
In this subsection we introduce some invariants of nanophrases over
$\alpha$
(some of them are defined in \cite{fu}). \par

Let $\Pi$ be the group which is defined as follows:
$$\Pi := (\{z_a\}_{a \in \alpha}|z_az_{\tau(a)}=1 \ for \ all \ a \in
\alpha).$$

\begin{definition}(cf. \cite{fu}).
Let $P=(\mathcal{A},(w_1|w_2|\cdots|w_k))$ be a nanophrase of
 length $k$ over $\alpha$ and $n_i$ the length of nanoword $w_i$. Set
$n=\sum_{1 \le i \le k}n_i$.
Then we define $n$ elements $\gamma_1^i$, $\gamma_2^i$, $\cdots$ and
 $\gamma_{n_i}^i$ ($i \in \{1,2,\cdots,k\}$) of $\Pi$ by
$\gamma_i^j := z_{|w_j(i)|}$ if $w_j(i) \neq w_l(m)$ 
 for  all  $l < j$   and   
	       for  all  $m < i$  when  $l=j$. 
 Otherwise $\gamma_i^j := z_{\tau(|w_j(i)|)}$.
Then we define $\gamma(P) \in \Pi^k$ by
$$\gamma(P):=(\gamma_1^1 \gamma_2^1 \cdots \gamma_{n_1}^1 , 
\gamma_1^2\gamma_2^2 \cdots \gamma_{n_2}^2 , \cdots ,
 \gamma_1^k \gamma_2^k \cdots \gamma_{n_k}^k).$$ 
\end{definition}

Then we obtain following proposition.

\begin{proposition}
$\gamma$ is a homotopy invariant of nanophrases.
\end{proposition}
   
We define a invariant of nanophrases $T$.\par
First we prepare some notations. Since the set $\alpha$ is a finite set,
we obtain following orbit decomposition of the $\tau$ : 
$\alpha/\tau = \{ \widetilde{a_{i_1}}, \widetilde{a_{i_2}}, \cdots 
,\widetilde{a_{i_l}}, \widetilde{a_{i_{l+1}}},\cdots, 
\widetilde{a_{i_{l+m}}} \}$, where 
$\widetilde{a_{i_j}}:= \{ a_{i_j}, \tau(a_{i_j}) \}$ such that 
$Card(\widetilde{a_{i_j}})=2$ for all $j \in \{1,\cdots,l \}$ and 
$Card(\widetilde{a_{i_j}})=1$ for all $j \in \{l+1,\cdots,l+m \}$ (we fix
a complete representative system $\{ a_{i_1}, a_{i_2}, \cdots 
,a_{i_l}, a_{i_{l+1}},\cdots, a_{i_{l+m}} \}$ which satisfy the above
condition). Let $\mathcal{A}$ be a $\alpha$-alphabet.  
For $A \in \mathcal{A}$ we define $\varepsilon(A) \in \{ \pm 1 \} $ by
$$\varepsilon(A):=\begin{cases} 1 \ ( \ if \ |A|=a_{i_j} \ for \ some \
		  j \in \{1,\cdots l+m \} \  ),
 \\ -1 \ ( \ if \ |A|=\tau({a_{i_j}}) \ for \ some \ j \in \{1,\cdots
		   l\} \ ). 
\end{cases}$$
Let $P=(\mathcal{A},(w_1|\cdots|w_k))$ be a nanophrase
over $\alpha$ and $A$, $B \in \mathcal{A}$. 
Let $K_{(i,j)}$ be $\mathbb{Z}$ if $i \le l$ and $j \le l$, 
otherwise $\mathbb{Z} / 2\mathbb{Z}$. 
We denote $K_{(1,1)} \times K_{(1,2)} \times \cdots K_{(1,l+m)} \times K_{(2,1)}
\times \cdots \times K_{(l+m,l+m)}$ by $\prod K_{(i,j)}$. 
Then we define
$\sigma_P(A,B) \in \prod {K_{(i,j)}}$ 
as follows:  
If $A$ and $B$ form
$\cdots A \cdots B \cdots A \cdots B \cdots$ in $P$, 
$|A| \in \widetilde{a_{i_{p}}}$and $|B|=a_{i_q}$
for some $m,n \in \{1, \cdots l+m \}$, or 
$\cdots B \cdots A \cdots B \cdots A \cdots$ in $P$, 
$|A| \in \widetilde{a_{i_{p}}}$ and
$|B|=\tau(a_{i_q})$ for some $p,q \in \{1, \cdots l+m\}$ , then
$\sigma_P(A,B):= (0,\cdots,0,\stackrel{(p,q)}{\check{1}},0,\cdots,0)$. 
If $\cdots A \cdots B \cdots A \cdots B \cdots$ in
 $P$, $|A| \in \widetilde{a_{i_{p}}}$ and  $|B|=\tau(a_{i_q})$, 
or $\cdots B \cdots A \cdots B \cdots A \cdots$
 in $P$, $|A| \in \widetilde{a_{i_{p}}}$ and $|B|=a_{i_q}$, then 
$\sigma_P(A,B):=(0,\cdots,0,\stackrel{(p,q)}{\check{-1}},0,\cdots,0)$. 
Otherwise $\sigma_P(A,B):=(0,\cdots,0)$. Under the above preparation, we
define the invariant $T$ as follows.

\begin{definition}
Let $P=(\mathcal{A}, (w_1|w_2|\cdots|w_k))$ be a nanophrase of
length $k$ over $\alpha$. For $A \in \mathcal{A}$ such that
 there exist $i \in \{1,2,\cdots,k\}$ with $Card(w_i^{-1}(A))=2$, 
we define $T_P(A) \in \prod K_{(i,j)}$ by
$$T_P(A):= \varepsilon(A)\sum_{B \in \mathcal{A}}\sigma_P(A,B),$$
and $T_P(w_i) \in \prod K_{(i,j)}$ by
$$T_P(w_i):= \sum_{A \in \mathcal{A},\ Card(w_i^{-1}(A))=2 }T_P(A).$$
Then we define $T(P) \in (\prod K_{(i,j)})^k$ by 
$$T(P):=(T_P(w_1),T_P(w_2),\cdots,T_P(w_k)).$$
\end{definition}

\begin{proposition}\label{T}
$T$ is a invariant of nanophrases over $\alpha$.
\end{proposition} 

\begin{proof}
It is clear that isomorphism does not change the value of $T$.
Consider the first homotopy move
$$P_{1}:=(\mathcal{A}, (xAAy)) \longrightarrow P_{2}:=(\mathcal{A} \setminus
\{A\}, (xy))$$
where $x$ and $y$ are words on $\mathcal{A}$, 
possibly including "$|$" character.
Since $A$ and $X$ are unlacement in the phrase $P_{1}$ 
for all $X \in \mathcal{A}$, $A$ dose not contribute to $T(P_{1})$.
So the first homotopy move does not change the value of $T$.\par
Consider the second homotopy move
$$P_{1}:=(\mathcal{A}, (xAByBAz)) \longrightarrow 
(\mathcal{A} \setminus \{A,B\}, (xyz))$$
where $|A|=\tau(|B|)$, and $x$, $y$ and $z$ are words on $\mathcal{A}$
possibly including "$|$" character. 
Suppose $y$ does not include "$|$" character and $Card(\widetilde{|A|})=2$
(So $Card(\widetilde{|B|})$ is also two). 
Then $T_{P_{1}}(A) + T_{P_{2}}(B) = 0$ since 
\begin{eqnarray*}
T_{P_1}(A)&=& \varepsilon(A)\left(\sigma_{P_1}(A,B)+\sum_{X \in \mathcal{A} 
              \setminus \{B\}}\sigma_{P_1}(A,X)\right) \\
          &=& \varepsilon(A)\sum_{X \in \mathcal{A} 
              \setminus \{B\}}\sigma_{P_1}(A,X) \\
          &=& -\varepsilon(B)\sum_{X \in \mathcal{A} 
              \setminus \{A\}}\sigma_{P_1}(B,X) \\
          &=& -\varepsilon(B)\left(\sigma_{P_1}(B,A)+\sum_{X \in \mathcal{A} 
              \setminus \{A\}}\sigma_{P_1}(B,X)\right) \\
          &=& -T_{P_1}(B).
\end{eqnarray*}

Moreover for 
$X \in \mathcal{A} \setminus \{A,B\}$, 
$\cdots A \cdots X \cdots A \cdots X \cdots$ \\
 (respectively $\cdots X \cdots A \cdots X \cdots A \cdots$) in $P_{1}$ 
if and only if  $\cdots B \cdots X \cdots B \cdots X \cdots$ 
 (respectively $\cdots X \cdots B \cdots X \cdots B \cdots$) in $P_{1}$.
and $|A| = \tau(|B|)$
So $\sigma_{P_{1}}(X,A)+\sigma_{P_{1}}(X,B)=0$ for all $X \in \mathcal{A}$.
So 

\begin{eqnarray*}
T_{P_1}(X)&=&\varepsilon(X)\left(\sigma_{P_1}(X,A)+\sigma_{P_1}(X,B)+
              \sum_{D \in \mathcal{A}\setminus \{A,B\}}\sigma_{P_1}(X,D)
               \right) \\
          &=& \varepsilon(X)\sum_{D \in \mathcal{A}\setminus \{A,B\}}
               \sigma_{P_1}(X,D) \\ 
          &=& \varepsilon(X)\sum_{D \in \mathcal{A}\setminus \{A,B\}}
               \sigma_{P_2}(X,D) \\ 
          &=& T_{P_2}(X). 
\end{eqnarray*}

This implies $T(P_{1}) = T(P_{2})$.\par
Suppose $y$ does not include "$|$" character and $Card(\widetilde{|A|})=1$
(So $Card(\widetilde{|B|})$ is also one).
This case also  $T_{P_{1}}(A) + T_{P_{2}}(B) = 0$ since 
\begin{eqnarray*}
T_{P_1}(A)&=& \varepsilon(A)\left(\sigma_{P_1}(A,B)+\sum_{X \in \mathcal{A} 
              \setminus \{B\}}\sigma_{P_1}(A,X)\right) \\
          &=& \varepsilon(A)\sum_{X \in \mathcal{A} 
              \setminus \{B\}}\sigma_{P_1}(A,X) \\
          &=& \varepsilon(B)\sum_{X \in \mathcal{A} 
              \setminus \{A\}}\sigma_{P_1}(B,X) \\
          &=& \varepsilon(B)\left(\sigma_{P_1}(B,A)+\sum_{X \in \mathcal{A} 
              \setminus \{A\}}\sigma_{P_1}(B,X)\right) \\
          &=& T_{P_1}(B),
\end{eqnarray*}
and all entry of $T_{P_{1}}(A)$ and $T_{P_{2}}(B)$ are elements of
$\mathbb{Z} / 2\mathbb{Z}$.
 Moreover for 
$X \in \mathcal{A} \setminus \{A,B\}$, $\cdots A \cdots X \cdots A \cdots X 
\cdots$ 
 (respectively $\cdots X \cdots A \cdots X \cdots A \cdots$) in $P_{1}$ 
if and only if  $\cdots B \cdots X \cdots B \cdots X \cdots$ 
 (respectively $\cdots X \cdots B \cdots X \cdots B \cdots$) in $P_{1}$.
Since $\widetilde{|A|}=\widetilde{|B|}$ and $Card(\widetilde{|A|})=1$
so 
$\sigma_{{P_{1}}}(X,A) = \sigma_{{P_{1}}}(X,B)$ in $\mathbb{Z} / 2 \mathbb{Z}$. 
So $\sigma_{P_{1}}(X,A)+\sigma_{P_{1}}(X,B)=0$ for all $X \in \mathcal{A}$.
By the above 
\begin{eqnarray*}
T_{P_1}(X)&=&\varepsilon(X)\left(\sigma_{P_1}(X,A)+\sigma_{P_1}(X,B)+
              \sum_{D \in \mathcal{A}\setminus \{A,B\}}\sigma_{P_1}(X,D)
               \right) \\
          &=& \varepsilon(X)
              \sum_{D \in \mathcal{A}\setminus \{A,B\}}\sigma_{P_1}(X,D) \\
          &=& \varepsilon(X)\sum_{D \in \mathcal{A}\setminus \{A,B\}}
               \sigma_{P_2}(X,D) \\ 
          &=& T_{P_2}(X). 
\end{eqnarray*}
This implies $T(P_{1}) = T(P_{2})$.\par
The case $y$ include "$|$" character is proved similarly.\par
Consider the third homotopy move 
$$P_{1}:= (\mathcal{A},(xAByACzBCt)) \rightarrow 
P_{2}:=(\mathcal{A},(xBAyCAzCBt))$$
where $|A|=|B|=|C|$, and $x$, $y$, $z$ and $t$ are words on 
$\mathcal{A}$ possibly including "$|$" character.
Suppose $y$ and $z$ do not including "$|$" character.
Note that $\sigma_{P_1}(A,B)=\sigma_{P_2}(A,C)$.
So
\begin{eqnarray*}
T_{P_1}(A)&=&\varepsilon(A)\left(\sigma_{P_1}(A,B)+\sum_{X \in \mathcal{A}
 \setminus \{B\}}\sigma_{P_1}(A,X)\right)\\
&=&\varepsilon(A)\left(\sum_{X \in \mathcal{A}
 \setminus \{C\}}\sigma_{P_2}(A,X)+\sigma_{P_2}(A,C)\right)\\
&=&T_{P_2}(A),
\end{eqnarray*}
and since $\sigma_{P_1}(C,B)=\sigma_{P_2}(C,A)$, we obtain 
\begin{eqnarray*}
T_{P_1}(C)&=&\varepsilon(C)\left(\sigma_{P_1}(C,B)+\sum_{X \in \mathcal{A}
 \setminus \{B\}}\sigma_{P_1}(C,X)\right)\\
&=&\varepsilon(C)\left(\sum_{X \in \mathcal{A}
 \setminus \{C\}}\sigma_{P_2}(C,X)+\sigma_{P_2}(C,A)\right)\\
&=&T_{P_2}(C).
\end{eqnarray*}
Moreover $\sigma_{P_1}(B,A)+\sigma_{P_1}(B,C)=0$ and 
$\sigma_{P_2}(B,A)=\sigma_{P_2}(B,C)=0$.
We obtain $T_{P_1}(B)=T_{P_2}(B)$.
It is checked easily that $T_{P_1}(E)=T_{P_2}(E)$ for all $E \neq A,B,C$.
So we obtain $T(P_1)=T(P_2)$.\par
The case $y$ or $z$ including "$|$" character is proved similarly.
\end{proof}

\begin{remark}
This invariant $T$ is the generalization of
invariants $T$ of nanophrases over $\alpha_0$ and the one-element set 
defined in \cite{fu}. If we use the invariant $T$ defined in this paper,
 then we can classify nanophrases of length 2 with 4 letters without 
the Lemma 4.2 in \cite{fu}.
\end{remark}

Next we define another new invariant.
Let $\pi$ be the group which is defined as follows:
$$\pi := (a \in \alpha|a\tau(a)=1 , ab = ba \ for \ all 
\ a, b \in \alpha \ ) \simeq \Pi / [\Pi,\Pi]. $$ 
Let $P=(\mathcal{A}, (w_1|w_2|\cdots|w_k))$ be a nanophrase of
length $k$ over $\alpha$. We define $(w_i , w_j)_P \in \pi$ for $i < j$
by
$$(w_i,w_j)_P := \prod_
{A \in Im(w_i) \cap Im(w_j)} |A|.$$  

\begin{proposition}\label{propair}
If nanophrases over $\alpha$, $P_1$ and $P_2$ are homotopic, then 
$(w_i,w_j)_{P_1} = (w_i,w_j)_{P_2}$. 
\end{proposition}
\begin{proof}
It is clear that isomorphisms does not change the value of $(w_i,w_j)_{P}$.
Consider the first homotopy move
$$P_{1}:=(\mathcal{A}, (xAAy)) \longrightarrow P_{2}:=(\mathcal{A} \setminus
\{A\}, (xy)).$$
In this move, the letter $A$ appear twice in the same component.
So $A$ dose not contribute to $(w_i,w_j)_{P_{1}}$.
This implies $(w_i,w_j)_{P_{1}} =  (w_i,w_j)_{P_{2}}$.\par
Consider second homotopy move
$$P_{1}:=(\mathcal{A}, (xAByBAz)) \longrightarrow 
(\mathcal{A} \setminus \{A,B\}, (xyz))$$
where $|A|=\tau(|B|)$, and $x$, $y$ and $z$ are words on $\mathcal{A}$
possibly including "$|$" character. 
Suppose $y$ does not include "$|$" character.
In this case, $A$ and $B$ are appear in the same component of nanophrase
$P_{1}$.
So $A$ and $B$ do not contribute to  $(w_i,w_j)_{P_{1}}$.
This implies $(w_i,w_j)_{P_{1}} =  (w_i,w_j)_{P_{2}}$ for all $i$, $j$.
Suppose $y$ include "$|$" character.
Suppose $A$ and $B$ are appear in the m-th component and the n-th component
of $P_{1}$.
Then 
\begin{eqnarray*}
(w_{m},w_{n})_{P_{1}} &=& (w_{m},w_{n})_{P_{2}} \cdot |A| \cdot |B| \\
                     &=& (w_{m},w_{n})_{P_{2}} \cdot |A| \cdot \tau(|A|)\\
                    &=& (w_{m},w_{n})_{P_{2}}, 
\end{eqnarray*}
and it is clear that $(w_{i},w_{j})_{P_{1}}=(w_{i},w_{j})_{P_{2}}$ 
for $(i,j) \neq (m,n)$.
So  $(w_{i},w_{j})_{P_{1}}=(w_{i},w_{j})_{P_{2}}$ for all $i$ and $j$. \par
Consider the third homotopy move 
$$P_{1}:= (\mathcal{A},(xAByACzBCt)) \rightarrow 
P_{2}:=(\mathcal{A},(xBAyCAzCBt))$$
where $|A|=|B|=|C|$, and $x$, $y$, $z$ and $t$ are words on 
$\mathcal{A}$ possibly including "$|$" character.
Note that the third homotopy move sent a letter in the l-th component
of $P_{1}$ to the l-th component of $P_{2}$.
So $(w_{i},w_{j})_{P_{1}}$ is not changed by the third homotopy move.\par
By the above, $(w_{i},w_{j})_{P_{1}}$ is a homotopy invariant of 
nanophrases.
\end{proof}
By the above proposition,
we obtain a homotopy invariant of nanophrases 
$$((w_{1},w_{2})_{P},(w_{1},w_{3})_{P},
\cdots,(w_{1},w_{k})_{P},(w_{2},w_{3})_{P},\cdots,(w_{k-1},w_{k})_{P}) 
\in \pi^{\frac{1}{2}k(k-1)}.$$

\subsection{The case of nanophrases of length more than or equal to $3$.}
Now using the invariants prepared in the last section and some lemmas, 
we classify the nanophrases of length more than or equal to $3$ with 
less than or equal to $4$ letters.
First recall the following lemmas from \cite{fu}.

\begin{lemma}
Let $P_1=(w_1|w_2|\cdots|w_k)$ and  $P_2=(v_1|v_2|\cdots|v_k)$ be 
nanophrases of length $k$ over $\alpha$. If $P_1$ and $P_2$ are
homotopic as nanophrases, then $w_i$ and $v_i$ are homotopic as 
\'etale words for all $i \in \{1,2,,\cdots,k\}$.
\end{lemma}

\begin{lemma}\label{mod2}
Let $P_1=(w_1|\cdots|w_k)$ and $P_2=(v_1|\cdots|v_k)$ 
be nanophrases of length $k$ over $\alpha$. 
If $P_1$ and $P_2$ are homotopic, then the length of $w_i$ is equal to 
length of $v_i$ modulo $2$ for all $i \in \{1,2,\cdots,k\}$.   
\end{lemma}

 A following lemma is checked easily by definition of homotopy of 
nanophrases.

\begin{lemma}\label{newlem}
Let $P_1=(w_1|\cdots|w_k)$ and $P_2=(v_1|\cdots|v_k)$ 
be nanophrases over $\alpha$.
If $P_1$ and $P_2$ are homotopic, then
 $(w_1|\cdots|w_{l}w_{l+1}|\cdots|w_k)$ and 
$(v_1|\cdots|v_{l}v_{l+1}|\cdots|v_k)$
are homotopic as nanophrases of length $k-1$ over $\alpha$ for all $l \in
 \{1,\cdots,k-1 \}$.  
\end{lemma} 

Now we give the classification theorem of nanophrases with $2$
letters. Set $P^{1,1;p,q}_{a} := 
(\emptyset|\cdots|\emptyset|\stackrel{p}{\check{A}}|\emptyset|
\cdots|\emptyset|\stackrel{q}{\check{A}}|\emptyset|\cdots|\emptyset)$
with $|A|=a$ for $1 \le p < q \le k$

\begin{theorem}
Let $P$ be a nanophrase of length $k$ with $2$ letters. 
Then $P$ is either homotopic to $(\emptyset|\cdots|\emptyset)$ 
or isomorphic to  $P^{1,1;p,q}_{a}$ for some $p,q \in \{1,\cdots
 k\}$, $a \in \alpha$. 
Moreover $P^{1,1;p,q}_{a}$ and 
$P^{1,1;p^ \prime,q^\prime}_{a^\prime}$ are homotopic 
if and only if $p=p^\prime$, $q=q^\prime$ and $a=a^\prime$.  
\end{theorem} 

\begin{proof}
The first part of this theorem is clear. 
We show the second part of this theorem.
By the definition of $(w_i,w_j)_P$, 
$(w_i,w_j)_{P^{1,1;p,q}_a} = a$ if $i=p$ and $j=q$. 
Otherwise $(w_i,w_j)_{P^{1,1;p,q}_a} = 1$.
For $a \in \alpha$, $a \neq 1$ in $\pi$. 
So if $P^{1,1;p,q}_a$ and $P^{1,1;p^\prime,q^\prime}_{a^\prime}$
are homotopic, then $p=p^\prime$, $q=q^\prime$ and $a=a^\prime$.   
\end{proof}

To describe the classification theorem of nanophrases with $4$
letters, we prepare following notations.\\
$P^{4;p}_{a,b}:=(\emptyset|\cdots|\emptyset|\stackrel{p}{\check{ABAB}
}|\emptyset|\cdots|\emptyset)$,\\
$P^{3,1;p,q}_{a,b}:=(\emptyset|\cdots|\emptyset|\stackrel{p}{\check{ABA}
}|\emptyset|\cdots|\emptyset|\stackrel{q}{\check{B}}|\emptyset|\cdots|
\emptyset)$,\\  
$P^{2,2I;p,q}_{a,b}:=(\emptyset|\cdots|\emptyset|\stackrel{p}{\check{AB}
}|\emptyset|\cdots|\emptyset|\stackrel{q}{\check{AB}}|\emptyset|\cdots|
\emptyset)$,\\
$P^{2,2II;p,q}_{a,b}:=(\emptyset|\cdots|\emptyset|\stackrel{p}{\check{AB}
}|\emptyset|\cdots|\emptyset|\stackrel{q}{\check{BA}}|\emptyset|\cdots|
\emptyset)$,\\
$P^{1,3;p,q}_{a,b}:=(\emptyset|\cdots|\emptyset|\stackrel{p}{\check{A}
}|\emptyset|\cdots|\emptyset|\stackrel{q}{\check{BAB}}|\emptyset|\cdots|
\emptyset)$,\\
$P^{2,1,1I;p,q,r}_{a,b}:=(\emptyset|\cdots|\emptyset|
\stackrel{p}{\check{AB}}|\emptyset|\cdots|\emptyset|\stackrel{q}{\check{A}}
|\emptyset|\cdots|\emptyset|\stackrel{r}{\check{B}
}|\emptyset|\cdots|\emptyset)$,\\
$P^{2,1,1II;p,q,r}_{a,b}:=(\emptyset|\cdots|\emptyset|
\stackrel{p}{\check{BA}}|\emptyset|\cdots|\emptyset|\stackrel{q}{\check{A}}
|\emptyset|\cdots|\emptyset|\stackrel{r}{\check{B}
}|\emptyset|\cdots|\emptyset)$,\\
$P^{1,2,1I;p,q,r}_{a,b}:=(\emptyset|\cdots|\emptyset|
\stackrel{p}{\check{A}}|\emptyset|\cdots|\emptyset|\stackrel{q}{\check{AB}}
|\emptyset|\cdots|\emptyset|\stackrel{r}{\check{B}
}|\emptyset|\cdots|\emptyset)$,\\
$P^{1,2,1II;p,q,r}_{a,b}:=(\emptyset|\cdots|\emptyset|
\stackrel{p}{\check{A}}|\emptyset|\cdots|\emptyset|\stackrel{q}{\check{BA}}
|\emptyset|\cdots|\emptyset|\stackrel{r}{\check{B}
}|\emptyset|\cdots|\emptyset)$,\\
$P^{1,1,2I;p,q,r}_{a,b}:=(\emptyset|\cdots|\emptyset|
\stackrel{p}{\check{A}}|\emptyset|\cdots|\emptyset|\stackrel{q}{\check{B}}
|\emptyset|\cdots|\emptyset|\stackrel{r}{\check{AB}
}|\emptyset|\cdots|\emptyset)$,\\
$P^{1,1,2II;p,q,r}_{a,b}:=(\emptyset|\cdots|\emptyset|
\stackrel{p}{\check{A}}|\emptyset|\cdots|\emptyset|\stackrel{q}{\check{B}}
|\emptyset|\cdots|\emptyset|\stackrel{r}{\check{BA}
}|\emptyset|\cdots|\emptyset)$,\\
$P^{1,1,1,1I;p,q,r,s}_{a,b}:=(\emptyset|\cdots|\emptyset|
\stackrel{p}{\check{A}}|\emptyset|\cdots|\emptyset|\stackrel{q}{\check{A}}
|\emptyset|\cdots|\emptyset|\stackrel{r}{\check{B}
}|\emptyset|\cdots|\emptyset|\stackrel{s}{\check{B}}
|\emptyset|\cdots|\emptyset)$,\\
$P^{1,1,1,1II;p,q,r,s}_{a,b}:=(\emptyset|\cdots|\emptyset|
\stackrel{p}{\check{A}}|\emptyset|\cdots|\emptyset|\stackrel{q}{\check{B}}
|\emptyset|\cdots|\emptyset|\stackrel{r}{\check{A}
}|\emptyset|\cdots|\emptyset|\stackrel{s}{\check{B}}
|\emptyset|\cdots|\emptyset)$,\\
$P^{1,1,1,1III;p,q,r,s}_{a,b}:=(\emptyset|\cdots|\emptyset|
\stackrel{p}{\check{A}}|\emptyset|\cdots|\emptyset|\stackrel{q}{\check{B}}
|\emptyset|\cdots|\emptyset|\stackrel{r}{\check{B}
}|\emptyset|\cdots|\emptyset|\stackrel{s}{\check{A}}
|\emptyset|\cdots|\emptyset)$,\\
with $|A|=a$, $|B|=b$. If $a=\tau(b)$, then nanophrases $P^{4;p}_{a,b}$,
$P^{2,2I;p,q}_{a,b}$ and $P^{2,2II;p,q}_{a,b}$ are 
homotopic to $(\emptyset|\cdots|\emptyset)$.
So when we write $P^{4;p}_{a,b}$,
$P^{2,2I;p,q}_{a,b}$, $P^{2,2II;p,q}_{a,b}$ we always assume 
that $a \neq \tau(b)$.\par   

Under the above notations the classification of nanophrases
 with $4$ letter is described as follows. 

\begin{theorem}\label{ctk4}
Let $P$ be a nanophrase of length $k$ with $4$ letters.
Then $P$ is either homotopic to nanophrase with less than or equal to $2$
 letters 
or isomorphic to $P^{X;Y}_{a,b}$ for some 
$X \in \{4,(3,1),\cdots,(1,1,1,1III)\}$, 
$Y \in \{1,\cdots,k,(1,2),\cdots,(k-3,k-2,k-1,k)\}$.
Moreover $P^{X;Y}_{a,b}$ and
 $P^{X^\prime;Y^\prime}_{a^\prime,b^\prime}$ are homotopic
if and only if $X=X^\prime$, $Y=Y^\prime$, $a=a^\prime$ and $b=b^\prime$.
\end{theorem}

\begin{proof}
The first part of this theorem is clear. 
We prove the rest of this theorem.
To prove Theorem \ref{ctk4}, it must be shown that (i) if $X \neq
 X^\prime$, then $P^{X;Y}_{a,b}$ and $P^{X^\prime;Y^\prime}_{a,b}$ are
 not homotopic; and (ii) each  of four letter nanophrase $P^{X;Y}_{a,b}$ 
is homotopic to $P^{X;Y^\prime}_{a,b}$ if and only if 
$Y=Y^\prime$, $a=a^\prime$ and $b=b^\prime$.
First we split basic shapes of nanophrases into 8 sets: 
$\mathcal{P}_0=\{(\emptyset|\cdots|\emptyset),P^{1,1;p}_{a}\}$,\\
$\mathcal{P}_1=\{P^{4;p}_{a,b}|1 \le p \le k , a,b \in \alpha \}$,\\ 
$\mathcal{P}_2=\{P^{3,1;p,q}_{a,b},P^{1,3;p,q}_{a,b}|1 \le p < q \le
 k, a,b \in \alpha\}$,\\ 
$\mathcal{P}_3=\{P^{2,2I;p,q}_{a,b},P^{2,2II;p,q}_{a,b}|1 \le p < q \le
 k, a,b \in \alpha\}$,\\
$\mathcal{P}_4=\{P^{2,1,1I;p,q,r}_{a,b}, P^{2,1,1II;p,q,r}_{a,b}|
1 \le p < q < r \le k , a,b \in \alpha\}$\\
$\mathcal{P}_5=\{P^{1,2,1I;p,q,r}_{a,b}, P^{1,2,1II;p,q,r}_{a,b}|
1\le p < q < r \le k , a,b \in \alpha\}$,\\
$\mathcal{P}_5=\{P^{1,1,2I;p,q,r}_{a,b}, P^{1,1,2II;p,q,r}_{a,b}|
1 \le p < q < r \le k , a,b \in \alpha\}$,\\
$\mathcal{P}_7=\{P^{1,1,1,1I;p,q,r,s}_{a,b}, 
P^{1,1,1,1II;p,q,r,s}_{a,b},P^{1,1,1,1III;p,q,r,s}_{a,b} |
1 \le p < q < r < s \le k , a,b \in \alpha\}$.\\
By using the invariants $\gamma$, $T$ and $((w_i,w_j)_P)_{i<j}$,  
we can easily check that 
two nanophrases $P \in \mathcal{P}_i$ and $P^\prime  \in \mathcal{P}_j$
 are homotopic only if $i=j$. This cuts down the number 
of pairs of nanophrases that need to be considered in (i).\par
Consider the nanophrases in $\mathcal{P}_1$.\par
The claim $P^{4;p}_{a,b}$ is homotopic to 
$P^{4;p^\prime}_{a^\prime,b^\prime}$ if and only if 
$p=p^\prime$, $a=a^\prime$ and $b=b^\prime$ follows from
Theorem \ref{tu} and Lemma \ref{newlem}.   
Consider the nanophrases in $\mathcal{P}_2$. \par
The claim $P^{3,1;p,q}_{a,b}$ is not homotopic to 
$P^{1,3;p^\prime,q^\prime}_{a^\prime,b^\prime}$: 
Suppose $P^{3,1;p,q}_{a,b}$ is homotopic to 
$P^{1,3;p^\prime,q^\prime}_{a^\prime,b^\prime}$.
Then $p = p^\prime$ and $q = q^\prime$, 
since
 $((w_i,w_j)_{P^{3,1;p,q}_{a,b}})_{i<j} =
((w_i,w_j)_{P^{1,3;p^\prime,q^\prime}_{a^\prime,b^\prime}})_{i<j}$.
By Lemma \ref{newlem} $(ABA|B)$ with $|A|=a, |B|=b$ 
must be homotopic to $(A^\prime |B^\prime
 A^\prime B^\prime)$
with $|A^\prime |=a^\prime, |B^\prime |=b^\prime$. 
However this contradicts Theorem \ref{ct24}.\par
The claim $P^{3,1;p,q}_{a,b}$ is homotopic to 
$P^{3,1;p^\prime,q^\prime}_{a^\prime,b^\prime}$ 
if and only if $p = p^\prime$, $q = q^\prime$, 
$a=a^\prime$ and $b = b^\prime$ follows by comparing
$((w_i,w_j)_{P^{3,1;p,q}_{a,b}})_{i<j}$ and 
$((w_i,w_j)_{P^{3,1;p^\prime,q^\prime}_{a^\prime,b^\prime}})_{i<j}$.\par
The claim $P^{1,3;p,q}_{a,b}$ is homotopic to 
$P^{1,3;p^\prime,q^\prime}_{a^\prime,b^\prime}$ 
if and only if $p = p^\prime$, $q = q^\prime$, 
$a=a^\prime$ and $b = b^\prime$ is proved similarly.\par     
Consider the nanophrases in $\mathcal{P}_3$. \par
The calm $P^{2,2I;p,q}_{a,b}$ and 
$P^{2,2II;p^\prime,q^\prime}_{a^\prime,b^\prime}$
are not homotopic:
Suppose $P^{2,2I;p,q}_{a,b}$ is homotopic to 
$P^{2,2II;p^\prime,q^\prime}_{a^\prime,b^\prime}$.
Then $p = p^\prime$ and $q = q^\prime$, 
since
 $((w_i,w_j)_{P^{2,2I;p,q}_{a,b}})_{i<j} =
((w_i,w_j)_{P^{2,2II;p^\prime,q^\prime}_{a^\prime,b^\prime}})_{i<j}$.      
By Lemma \ref{newlem} $(AB|AB)$ with $|A|=a, |B|=b$ 
must be homotopic to $(A^\prime B^\prime |B^\prime A^\prime)$
with $| A^\prime |=a^\prime, | B^\prime |=b^\prime$. 
However this contradicts Theorem \ref{ct24}.\par
The claim $P^{2,2I;p,q}_{a,b}$ and 
$P^{2,2I;p^\prime,q^\prime}_{a^\prime,b^\prime}$
are homotopic if and only if $p = p^\prime$, $q = q^\prime$ 
$a=a^\prime$ and $b = b^\prime$ follows by comparing values of the 
invariant $((w_i,w_j)_P)_{i<j}$.\par
The claim $P^{2,2II;p,q}_{a,b}$ and 
$P^{2,2II;p^\prime,q^\prime}_{a^\prime,b^\prime}$
are homotopic if and only if $p = p^\prime$, $q = q^\prime$ 
$a=a^\prime$ and $b = b^\prime$ is proved similarly.\par
Consider the nanophrases in $\mathcal{P}_4$.\par
The claim $P^{2,1,1I;p,q,r}_{a,b}$ and 
$P^{2,1,1II;p^\prime,q^\prime,r^\prime}_{a^\prime,b^\prime}$
are not homotopic:
Suppose $P^{2,1,1I;p,q,r}_{a,b}$ is homotopic to 
$P^{2,1,1II;p^\prime,q^\prime,r^\prime}_{a^\prime,b^\prime}$.
Then $p = p^\prime$, $q = q^\prime$ and $r=r^\prime$ 
since
 $((w_i,w_j)_{P^{2,1,1I;p,q}_{a,b}})_{i<j} =
((w_i,w_j)_{P^{2,1,1II;p^\prime,q^\prime,r^\prime}
_{a^\prime,b^\prime}})_{i<j}$.
By Lemma \ref{newlem} nanophrases $(ABA|B)$ and $(B^\prime
 A^\prime A^\prime | B^\prime)$ are homotopic.
However this contradicts Theorem \ref{ct24}.\par
The claim $P^{2,1,1I;p,q,r}_{a,b}$ and 
$P^{2,1,1I;p^\prime,q^\prime,r^\prime}_{a^\prime,b^\prime}$ are
homotopic if and only if  $p=p^{\prime}, q=q^{\prime}$ and
$r=r^{\prime}$,
 $a=a^\prime$ and $b=b^\prime$ follows by comparing values of the 
invariant $((w_i,w_j)_P)_{i<j}$.\par
For the nanophrases in $\mathcal{P}_5$ and $\mathcal{P}_6$, we can prove
(i) and (ii) similarly.\par
Consider the nanophrases in $\mathcal{P}_7$.\par 
The claim nanophrases $P^{1,1,1,1I;p,q,r,s}_{a,b}$ and  
$P^{1,1,1,1II;p^\prime,q^\prime,r^\prime,s^\prime}_{a^\prime,b^\prime}$
are not homotopic: 
Indeed if we assume $P^{1,1,1,1I;p,q,r,s}_{a,b}$ and 
$P^{1,1,1,1II;p^\prime,q^\prime,r^\prime,s^\prime}_
{a^\prime,b^\prime}$ are homotopic, then 
 $p=p^{\prime}, q=q^{\prime},
r=r^{\prime}$ and $z=z^{\prime}$ 
since $((w_i,w_j)_{P^{1,1,1,1I;p,q,r,s}_{a,b}})_{i<j}=
((w_i,w_j)_{P^{1,1,1,1II;p^\prime,q^\prime,r^\prime,s^\prime}_
{a^\prime,b^\prime}})_{i<j}$. So $(A|BAB)$ must be homotopic to
 $(A^\prime |A^\prime B^\prime B^\prime)$ by Lemma \ref{newlem}. 
But this contradicts Theorem \ref{ct24}.\par
The claim nanophrases $P^{1,1,1,1I;p,q,r,s}_{a,b}$ and  
$P^{1,1,1,1III;p^\prime,q^\prime,r^\prime,s^\prime}_{a^\prime,b^\prime}$
are not homotopic: 
 If we assume $P^{1,1,1,1I;p,q,r,s}_{a,b}$ and 
$P^{1,1,1,1III;p^\prime,q^\prime,r^\prime,s^\prime}_
{a^\prime,b^\prime}$ are homotopic, then $p=p^{\prime}, q=q^{\prime},
r=r^{\prime}$ and $z=z^{\prime}$, then $(A|AB|B)$ must be homotopic to
 $(A^\prime |\emptyset|A^\prime )$ by Lemma \ref{newlem}.
However this contradicts the homotopy invariance of
 $((w_i,w_j)_P)_{i<j}$.\par
The claim nanophrases $P^{1,1,1,1II;p,q,r,s}_{a,b}$ and  
$P^{1,1,1,1III;p^\prime,q^\prime,r^\prime,s^\prime}_{a^\prime,b^\prime}$
are not homotopic follows similarly as the above.\par
The claim nanophrases $P^{1,1,1,1I;p,q,r,s}_{a,b}$ and 
$P^{1,1,1,1I;p^\prime,q^\prime,r^\prime,s^\prime}_{a^\prime,b^\prime}$
are homotopic if and only if 
 $p=p^{\prime}, q=q^{\prime}, r=r^{\prime}$ and $z=z^{\prime}$
, $a=a^\prime$ and $b=b^\prime$ follows by homotopy invariance of 
$((w_i,w_j)_P)_{i<j}$.
The claim nanophrases $P^{1,1,1,1II;p,q,r,s}_{a,b}$ and 
$P^{1,1,1,1II;p^\prime,q^\prime,r^\prime,s^\prime}_{a^\prime,b^\prime}$
are homotopic if and only if 
 $p=p^{\prime}, q=q^{\prime},
r=r^{\prime}$ and $z=z^{\prime}$, $a=a^\prime$ and $b=b^\prime$ and 
the claim nanophrases $P^{1,1,1,1III;p,q,r,s}_{a,b}$ and 
$P^{1,1,1,1III;p^\prime,q^\prime,r^\prime,s^\prime}_{a^\prime,b^\prime}$
are homotopic if and only if 
 $p=p^{\prime}, q=q^{\prime},
r=r^{\prime}$ and $z=z^{\prime}$, 
$a=a^\prime$ and $b=b^\prime$ follows similarly.\par
Now the we have completed the homotopy classification of nanophrases
with less than of equal to four letters without the condition on length. 
\end{proof}

\section{Proof of The Theorem \ref{mthm}.}
To complete the proof of the Theorem \ref{mthm}, we prepare a following lemma.

\begin{lemma}\label{lem}
The nanophrases over $\alpha$, $(A|A)$, $(AB|AB)$ with  $|A| \neq \tau(|B|)$ , 
$(AB|BA)$ with $|A| \neq \tau(|B|)$, $(ABA|B)$, $(A|BAB)$, 
$(AB|A|B)$, 
$(BA|A|B)$, 
$(A|AB|B)$,
$(A|BA|B)$, $(A|B|AB)$, $(A|B|BA)$, $(A|A|B|B)$, $(A|B|A|B)$ and 
$(A|B|B|A)$ are not homotopic to nanophrases over $\alpha$ 
which have the empty words in its components.
\end{lemma}
\begin{proof}
This lemma easily follows from Proposition \ref{propair}, 
Lemma \ref{mod2} and Theorem \ref{ct24}.   
\end{proof}

Now Theorem \ref{mthm} immediately follows from 
Theorem \ref{ctk4} and Lemma \ref{lem}.
It is sufficient to apply the above theorems to the case
$\alpha=\alpha_0$ with involution $\tau:\alpha_0
 \rightarrow \alpha_0$ permuting $a$ and $b$.\par

\vspace{0.3cm}
Department of Mathematics, Hokkaido University

Sapporo 060-0810, Japan

email: fukunaga@math.sci.hokudai.ac.jp


\begin{thebibliography}{9}
\bibitem[1]{fu} T. Fukunaga, \emph{Homotopy classification of
	   nanophrases in Turaev's theory of words}, to appear in
	   Journal of Knot Theory and Its Ramifications (2009). 
\bibitem[2]{ga} C. F. Gauss, \emph{Werke},  Vol.8, Teubner, Leipzig,
	    1900.
\bibitem[3]{kad} T. Kadokami, Detecting non-triviality of virtual links,
           Journal of Knot Theory and Its Ramifications \textbf{12}
           (2003), no. 6, 781-803. 
\bibitem[4]{ka2} N.Kamada and S.Kamada, \emph{Abstract link diagrams and
	   virtual knots}, Journal of Knot Theory and Its Ramifications 
           \textbf{9} (2000), no.1, 93-106.
\bibitem[5]{SW} D. S. Silver, S. G. Williams, 
           \emph{An Invariant for open virtual strings}, 
           Journal of Knot Theory and Its Ramifications \textbf{15} (2006),
           no.2, 143-152.      
\bibitem[6]{tu1} V. Turaev, \emph{Topology of words}, Proceedings of
	   the London Mathematical Society \textbf{95} (2007), no.2, 360-417.
\bibitem[7]{tu2} V. Turaev, \emph{Knots and words},
	     International Mathematics Research Notices (2006), 
             Art. ID 84098, 23 pp. 
\bibitem[8]{tu3} V. Turaev, \emph{Lectures on topology of words},
	     Japanese Journal of Mathematics \textbf{2} (2007), 1-39.  
\bibitem[9]{tu5} V. Turaev, \emph{Virtual strings}, Annals de
	     l'Institut Fourier \textbf{54} (2004), no.7, 2455-2525.  
\end{thebibliography}
\end{document}